\documentclass{svcyclop}

\usepackage{graphicx}    % standard LaTeX graphics tool
\usepackage{amsmath, amsfonts, amssymb}
\newcommand{\be}{\begin{equation}}
\newcommand{\ee}{\end{equation}}
\newcommand{\bea} {\begin{eqnarray}}
\newcommand{\nbea} {\begin{eqnarray*}}
\newcommand{\eea} {\end{eqnarray}}
\newcommand{\neea} {\end{eqnarray*}}

\def\RR{{R\!\!\!\!\!I~}}

\def\del{\partial}

\def\mm{{\bf m}}
\def\xx{{\bf x}}

\def\vv{{\bf v}}

\def\aa{{\bf a}}
\def\bb{{\bf b}}
\def\ff{{\bf f}}

\def\xxi{{\bf \xi}}

\def \uu{{\bf u}}

\begin{document}

\title{Gas Dynamics Equations: Computation}

\author{Gui-Qiang G. Chen} 

\institute{Mathematical Institute, University of Oxford\\24--29 St Giles, Oxford, OX1 3LB, United Kingdom\\
Homepage: http://people.maths.ox.ac.uk/chengq/ \\
Email: chengq@maths.ox.ac.uk
}

\maketitle

Shock waves, vorticity waves, and entropy waves are
fundamental discontinuity waves in nature
and arise in supersonic or transonic gas flow, or from a very sudden
release (explosion) of chemical, nuclear, electrical, radiation,
or mechanical energy in a limited space.
Tracking these discontinuities and their interactions,
especially when and where new waves arise and
interact in the motion of gases, is
one of the main motivations for numerical computation for the gas dynamics equations.

The fundamental equations governing the dynamics of gases are the compressible
Euler equations, consisting of conservation laws of mass, momentum, and
energy:
\begin{equation}\label{Euler1}
\del_t\rho +\nabla\cdot \mm=0,\,\,\,
\del_t\mm+\nabla\cdot\big(\frac{\mm\otimes\mm}{\rho}\big)+\nabla p=0,\,\,\,
\del_t (\rho E) +\nabla\cdot\big(\mm(E+p)\big)=0,
\end{equation}
where $\nabla$ is the gradient with respect to the space
variable $\xx\in\RR^d$,
$\rho$ is the density,
$\vv\in \RR^d$ is the gas velocity with
$\rho\vv=\mm$ the momentum vector, $p$ is the scalar pressure, and
$E=\frac{1}{2}|\vv|^2+ e(\tau,p)$ is the total energy
with $e$ the internal energy, a given function of $(\rho,p)$
defined through thermodynamical relations.
The notation $\aa\otimes\bb$ denotes the tensor product of two vectors.
The other two thermodynamic variables are the temperature $\theta$ and
the entropy $S$. If $(\rho, S)$ are chosen as the independent variables,
then the constitutive relations
$(e,p,\theta)=(e(\rho,S), p(\rho,S),\theta(\rho,S))$
are governed by
$\theta dS=de +p \,d(\frac{1}{\rho})$.
For a polytropic gas,
$p=R\rho\theta, e=c_v\theta,\gamma=1+\frac{R}{c_v}$,
and
\begin{equation}
p=p(\rho,S)=\kappa\rho^\gamma e^{S/c_v}, \qquad
e=\frac{\kappa}{\gamma-1}\rho^{\gamma-1}e^{S/c_v}
   =\frac{\theta}{\gamma-1},
\label{1.6}
\end{equation}
where $R, c_v$, and $\kappa$ are positive constants, respectively.
System \eqref{Euler1} is complemented by the Clausius inequality:
$\del_t(\rho a(S))+\nabla\cdot ({\mm} a(S))\ge 0$
in the sense of distributions for any $a(S)\in C^1, a'(S)\ge 0$,
to identify physical shocks.

The Euler equations for an isentropic gas take the
simpler form:
\begin{equation}
\del_t\rho +\nabla\cdot \mm=0,\quad
\del_t\mm+\nabla\cdot\big(\frac{\mm\otimes\mm}{\rho}\big)
+\nabla p=0,
\label{Euler2}
\end{equation}
where
$p(\rho)=\kappa_0\rho^\gamma$ with constants $\gamma>1$ and
$\kappa_0>0$.

These systems fit into the general form of hyperbolic conservation laws:
\begin{equation}
\del_t \uu +\nabla \cdot \ff(\uu)=0, \qquad \uu\in\RR^m,\,\, \xx\in \RR^d,
\label{cons}
\end{equation}
where $\ff: \RR^m\to (\RR^m)^d$ is a nonlinear mapping.
Besides \eqref{Euler1} and \eqref{Euler2}, most of partial differential
equations arising from physical or engineering science can be also
formulated into form \eqref{cons} or its variants, for example, with
additional source terms or equations modeling the effects of dissipation,
relaxation, memory, damping, dispersion, magnetization, etc.
Hyperbolicity of  system \eqref{cons} requires that,
for all $\xxi\in S^{d-1}$, the matrix $(\xxi\cdot\nabla \ff(\uu))_{m\times m}$
have $m$ real eigenvalues $\lambda_j(\uu,\xxi), j=1,2,\cdots,m,$ and be
diagonalizable.

The main difficulty in calculating fluid flows with discontinuities is that
it is very hard to predict, even in the process of a flow calculation,
when and where new discontinuities arise and interact. Moreover, tracking the discontinuities,
especially their interactions, is numerically burdensome
(see \cite{BCLW,Da,HR,La3}).

One of the efficient numerical approaches is shock capturing algorithms.
Modern numerical ideas of shock capturing for computational
fluid dynamics can date back to 1944 when
von Neumann first proposed a new numerical method,
a centered difference scheme,
to treat the hydrodynamical shock problem, for which
numerical calculations showed oscillations on mesh
scale (see Lax \cite{La1}).
von Neumann's dream of capturing shocks was first realized when
von Neumann and Richtmyer \cite{NR} in 1950 introduced
the ingenious idea of adding a numerical
viscous term of the same size as the truncation error
into the hydrodynamic equations.
Their numerical viscosity guarantees that the scheme is consistent with
the Clausius inequality, i.e., the entropy inequality.
The shock jump conditions, the Rankine-Hugoniot jump conditions,
are satisfied, provided that the Euler equations of gas dynamics are discretized
in conservation form.
Then oscillations were eliminated by the judicious use of
the artificial viscosity; solutions constructed by this method converge
uniformly except in a neighborhood of shocks,
where they remain bounded and are spread out over a few mesh intervals.

\medskip
Related analytical idea of shock capturing, vanishing viscosity methods,
is quite old.
For example, there are some hints about the idea of regarding inviscid
gases as viscous gases with vanishingly small viscosity in the seminal
paper by Stokes (1848), as well as the important contributions
of Rankine (1870),
Hugoniot (1989), and Rayleigh (1910).
See Dafermos \cite{Da} for the details.

\medskip
The main challenge in designing shock capturing numerical algorithms is
that weak solutions are not unique; and the numerical schemes should be
consistent with the Clausius inequality, the entropy inequality.
Excellent numerical schemes should be also numerically simple, robust, fast,
and low cost, and have sharp oscillation-free resolutions and high
accuracy in domains where the solution is smooth.
It is also desirable that the schemes capture vortex sheets, vorticity waves,
and entropy waves, and are coordinate invariant, among others.

For the one-dimensional case, examples of success
include the Lax-Friedrichs scheme (1954), the Glimm scheme (1965),
the Godunov scheme (1959) and related high order schemes; for example,
van Leer's MUSCL (1981), Colella-Wooward's PPM (1984),
Harten-Engquist-Osher-Chakravarthy's ENO (1987), the more recent WENO (1994, 1996),
and the Lax-Wendroff scheme (1960) and its two-step version, the Richtmyer
scheme (1967) and the MacCormick scheme (1969). See
\cite{CLiu,CToro,Da,FJ,GR,Le,TLZ,To} and the references cited therein.

For the multi-dimensional case, one direct approach is to generalize directly
the one-dimensional methods to solve multi-dimensional problems;
such an approach has
led several useful numerical methods including semi-discrete methods
and Strang's dimension-dimension splitting methods.

Observe that multi-dimensional effects do play a significant role in the
behavior of the solution locally, and the approach that only solves
one-dimensional Riemann problems in the coordinate directions
clearly lacks the use of all the multi-dimensional information.
The development of fully multi-dimensional methods requires
a good mathematical theory to understand the multi-dimensional behavior of
entropy solutions; current efforts in this direction
include using more information about the multi-dimensional behavior
of solutions, determining the direction of primary wave propagation and
employing wave propagation in other directions,
and using transport techniques, upwind techniques,
finite volume techniques, relaxation techniques, and kinetic
techniques from the microscopic level.
See \cite{CCY1,KT,LL2,TLZ}.
Also see \cite{FJ,GM,GR,Le,To} and the references cited therein.

Other useful methods to calculate sharp fronts for gas dynamics equations include
front-tracking algorithms \cite{CGM,GKM}, level set methods \cite{Osher,Sethian},
among others.

%%%%%%%%%%%%%%%%%%%%%%%%%%%%%%%%%%%%%%%%%%%%%%%%%%%%%%%%%%%%%%%%%%%%%%

\end{document}